\documentclass[12pt,twoside]{article}
\usepackage{amsmath,amsbsy,amsfonts}
\usepackage{color}
\newtheorem{theorem}{Theorem}[section]

\let\Section=\section
\def\section{\setcounter{equation}{0}\Section}

\begin{document}

\date{}
\title{Existence of solutions for a class of singular elliptic systems with
convection term}
\author{\textsf{{Claudianor O. Alves} \thanks{%
C.O. Alves was partially supported by CNPq/Brazil 303080/2009-4,
coalves@dme.ufcg.edu.br}} \\
{\small \textit{Unidade Acad\^emica de Matem\'atica e Estat\'{\i}stica}}\\
{\small \textit{Universidade Federal de Campina Grande}}\\
{\small \textit{58429-900, Campina Grande - PB - Brazil}}\\
{\small \textit{e-mail address: coalves@dme.ufcg.edu.br}} \\
\\
\vspace{1mm} \textsf{{Abdelkrim Moussaoui} \thanks{%
A. Moussaoui was supported by the European program Averro\`es-Erasmus
Mundus(Grant No $\sharp $1872)}}\\
{\small \textit{Biology Department, A. Mira Bejaia University, }}\\
{\small \textit{Targa Ouzemour 06000 Bejaia, Algeria}}\\
{\small \textit{e-mail address:remdz@yahoo.fr}}}
\maketitle

\begin{abstract}
We show the existence of positive solutions for a class of singular elliptic
systems with convection term. The approach combines sub
and supersolution method with the pseudomonotone operator theory and perturbation arguments involving singular terms.
\end{abstract}



\section{Introduction}

\label{S1}

In this work, we focus our attention on the existence of solutions for the
following class of elliptic system with convection term
\begin{equation*}
\left\{
\begin{array}{l}
-\Delta {u}=\frac{1}{v^{\alpha _{1}}}\pm v^{\beta _{1}}+g_{1}(\nabla
u,\nabla v)\,\,\,\mbox{in}\,\,\,\Omega , \\
\mbox{} \\
-\Delta {v}=\frac{1}{u^{\alpha _{2}}}\pm u^{\beta _{2}}+g_{2}(\nabla
u,\nabla v)\,\,\,\mbox{in}\,\,\,\Omega , \\
\mbox{} \\
u,v>0\,\,\ \mbox{in}\,\,\,\Omega , \\
\mbox{} \\
u=v=0\,\,\ \mbox{on}\,\,\, \partial  \Omega , \\
\end{array}%
\right. \leqno{(S)_\pm}
\end{equation*}%
where $\Omega $ is a bounded domain with smooth boundary and $g_{i}:\mathbb{R%
}^{2N}\rightarrow \lbrack 0,+\infty ),$ $i=1,2,$ are positive continuous
functions belong to $L^{\infty }(\mathbb{R}^{2N})$. We consider the system $%
(S)_{\pm }$ in a singular case assuming $\alpha _{i},\beta _{i},\in \lbrack
0,1)$ for $i=1,2.$

Hereafter $(u,v)$ is a solution to $(S)_{\pm }$ if $u,v\in C^{2}(\Omega
)\cap H_{0}^{1}(\Omega )$ are both positive in $\Omega $ and satisfy the
equations of $(S)_{\pm }$ in the classical sense.

Nonlinear singular boundary value problems are mathematically challenging
and important for applications. They arise in several physical situations
such as fluid mechanics, pseudoplastics flow, chemical heterogeneous
catalysts, non- Newtonian fluids, biological pattern formation, for more
details about this subject, we cite the papers of Fulks \& Maybe \cite{FM},
Callegari \& Nashman \cite{CN1,CN2} and the references therein.

Systems $(S)_{\pm }$ can be see as a version of the singular scalar
equations
\begin{equation*}
\left\{
\begin{array}{l}
-\Delta {u}=\frac{1}{u^{\alpha }}\pm u^{\beta }+g(\nabla u)\,\,\,\mbox{in}%
\,\,\,\Omega , \\
\mbox{} \\
u>0\,\,\ \mbox{in}\,\,\,\Omega , \\
\mbox{} \\
u=0\,\,\ \mbox{on} \,\,\, \partial   \Omega ,%
\end{array}%
\right. \leqno{(P)_\pm}
\end{equation*}%
with $\alpha ,\beta >0$ and $g:\mathbb{R}^{N}\rightarrow \mathbb{R}$ be a
continuous function verifying some technical conditions. Several works are
devoted to classes of problems covering $(P)_{\pm }$. For instance, see the
papers of Aranda \cite{AR}, Ghergu \& Radulescu \cite{GR1,GR2}, Giarrusso \&
Porru \cite{GP}, Lair \& Wood \cite{LW}, Zhang \cite{Z} and references
therein. Problem $(P)_{\pm }$ without a convection term, that is $g=0$ was
also investigated. Relevant contributions regarding this situation can be
found in Crandall, Rabinowitz \& Tartar \cite{CRT}, Choi \& McKenna \cite{CM}%
, Coclite \& Palmieri \cite{CP}, C\'{\i}rstea, Ghergu \& Radulescu \cite{CGR}%
, D\'{a}vila \& Montenegro \cite{DM} and Diaz, Morel \& Oswald \cite{DMO}.
The main tools used in the aforementioned works are Sub and Supersolution,
Fixed Point Theorems, Bifurcation Theory and Galerkin Method. On the other
hand, using variational technique, more precisely mountains pass theorem, de
Figueiredo, Girardi \& Matzeu \cite{FGM} studied a class of elliptic
problems without singularity, where the nonlinearity depends of the gradient
of the solution.

Related to systems $(S)_{\pm }$, to date, the only case considered in the
literature, known to authors for $g_{i}\not=0$ is the paper due to Alves,
Carri\~{a}o \& Faria \cite{ACL}. For the case where $g_{i}=0$, we refer the
reader to the survey paper by Alves \& Corr\^{e}a \cite{AlvesCorrea}, Alves,
Corr\^{e}a \& Gon\c{c}alves \cite{ACG}, El Manouni, Perera \& Shivaji \cite%
{EPS} and Motreanu \& Moussaoui \cite{MM}. From the above commentaries, we
observe that in recent years singular elliptic problems with convection term
has received few attention. Motivated by this fact, our aim is to show the
existence of solutions for a class of elliptic systems where the
nonlinearity besides a singular term has a convection term. The proof
combines results involving pseudomonotone operators, sub and supersolution method and perturbation
arguments involving singular terms. We emphasize that our study complete
those made in \cite{AlvesCorrea}, \cite{ACG} and \cite{MM}, in the sense
that in those papers the authors did not consider the case where the nonlinearity has a
convection term, and also \cite{ACL}, because a different type of singular
term was considered. The method used in the present work is different from
those applied in the aforementioned papers.

Our main result is the following:

\begin{theorem}
\label{T1} Assume that $g_{i}:\mathbb{R}^{N}\rightarrow \mathbb{R}$ are
continuous functions and $\alpha _{i},\beta _{i}\in \lbrack 0,1)$ for $i=1,2$%
. Then, the systems $(S)_{\pm }$ has a solution.
\end{theorem}

The proof of Theorem \ref{T1} is done in Sections \ref{S3} and \ref{S4}. The
first main technical difficulty is that the nonlinearities of $(S)_{\pm }$
depend of gradient of the solution, which is more a complicating factor.
Indeed, for the scalar case, an interesting result is proved by Kazdan \&
Kramer \cite{KK} and Leon \cite{leon}, where the authors develop a sub and supersolution method
for scalar problem where the nonlinearity depends of the gradient. Instead,
the counter-part of this result for systems with gradient terms is not known
in the literature. Thus, we do not know a result involving sub and
supersolution that could be use to establishes the existence of solution for
this class of system. To overcome this difficulty, we show in Section \ref%
{S2} a result, see Theorem \ref{T2}, which can be see as a sub and
supersolution method for systems whose nonlinearity depends of gradient.

The second main difficulty in the proof of Theorem \ref{T1} is associated
with the fact that the sub and supersolution method in its version involving
maximum principle cannot be used directly for systems involving the gradient
of the solution. Moreover, the way as the singularities appear in the system
$(S)_{\pm }$ is a difficult point to work with maximum principle. In order
to overtake the stated problem we first introduce a parameter $\varepsilon
>0 $ in $(S)_{\pm }$, giving rise to regularized systems for $(S)_{\pm }$\
whose study is relevant for our initial problem. Then, for the regularized
systems, we combine variational methods with the sub-supersolution one to
prove the existence of a solution $(u_{\varepsilon },v_{\varepsilon })\in
H^{1}(\Omega )\cap L^{\infty }(\Omega )\times H^{1}(\Omega )\cap L^{\infty
}(\Omega )$. This solution $(u_{\varepsilon },v_{\varepsilon })$ is located
in some rectangle formed by the sub and supersolution, independent for $%
\varepsilon >0$, which does not contains zero for all $\varepsilon >0$.
Then, a positive solution of $(S)_{\pm }$ is obtained by passing to the
limit as $\varepsilon \rightarrow 0$. This is based on a priori estimates
and Hardy-Sobolev inequality. The positivity of the solution is derived from
the independence of the subsolution of the regularized systems on $%
\varepsilon $.

The rest of this article is organized as follows: In section \ref{S2} we
state and prove a general theorem about sub and supersolution method for
systems with convection term. Sections \ref{S3} and \ref{S4} contain the
proof of Theorem \ref{T1}.

\section{ An auxiliary result}

\label{S2}

The main goal in this section is to prove the Theorem \ref{T2} below, which
is a key point in the proof of Theorem \ref{T1}. An interesting point
related to Theorem \ref{T2} is the fact that it is a result of sub and
supersolution whose the proof is made using pseudomonotone operator theory. 

\begin{theorem}
\label{T2} Let $H,G:\Omega \times \mathbb{R}^{+}\times \mathbb{R}^{+}\times \mathbb{R}%
^{N}\times \mathbb{R}^{N}\rightarrow \mathbb{R}$ continuous functions 
function verifying the following conditions: Given $T,S>0$, there exist $C>0$ and  $\alpha, \beta \in (0,1)$, such that
$$
|H(x,s,t,\eta,\xi)|,|G(x,s,t,\eta,\xi)| \leq C(1+|\eta|^{\alpha}+|\xi|^{\beta})   
$$
for all $(x,s,t,\eta,\xi) \in \Omega \times [0,T] \times [0,S] \times  \mathbb{R}^{N} \times \mathbb{R}^{N}$. Let $\tilde{g},\widehat{g}\in C^{2}(\overline{\Omega })$ and $\underline{u},\overline{u},\underline{v},\overline{%
v}\in W^{1,\infty}(\Omega ) $ with
\begin{equation*}
\underline{u}(x)\leq \tilde{g}(x)\leq \overline{u}(x)\,\,\,\mbox{on}%
\,\,\,\partial \Omega
\end{equation*}%
and
\begin{equation*}
\underline{v}(x)\leq \widehat{g}(x)\leq \overline{v}(x)\,\,\,\mbox{on}%
\,\,\,\partial \Omega .
\end{equation*}
Assume that 
\begin{equation*}
\int_{\Omega }\nabla \underline{u}\,\nabla \phi dx\leq \int_{\Omega }H(x,%
\underline{u},\underline{v},\nabla \underline{u},\nabla \underline{v})\phi \,dx,
\end{equation*}%
\begin{equation*}
\int_{\Omega }\nabla \underline{v}\,\nabla \psi dx\leq \int_{\Omega }G(x,%
\underline{u},\underline{v},\nabla \underline{u},\nabla \underline{v})\psi \,dx,
\end{equation*}%
\begin{equation*}
\int_{\Omega }\nabla \overline{u}\,\nabla \phi dx\geq \int_{\Omega }H(x,%
\overline{u},\overline{v},\nabla \overline{u},\nabla \overline{v})\phi \,dx
\end{equation*}%
and
\begin{equation*}
\int_{\Omega }\nabla \overline{v}\,\nabla \psi dx\geq \int_{\Omega }G(x,%
\overline{u},\overline{v},\nabla \overline{u},\nabla \overline{v})\psi \,dx,
\end{equation*}%
for all nonnegative functions $\phi ,\psi \in H^{1}(\Omega )$. Then, there is
\linebreak $(u,v)\in (H^{1}(\Omega )\cap L^{\infty }(\Omega ))\times
(H^{1}(\Omega )\cap L^{\infty }(\Omega ))$ verifying
\begin{equation*}
\underline{u}(x)\leq u(x)\leq \overline{u}(x)\,\,\,\mbox{and}\,\,\,%
\underline{v}(x)\leq v(x)\leq \overline{v}(x)\,\,\forall x\in \Omega ,
\end{equation*}%
\begin{equation*}
u-\tilde{g},v-\widehat{g}\in H_{0}^{1}(\Omega )
\end{equation*}%
and
\begin{equation*}
\begin{array}{l}
\displaystyle\int_{\Omega }\nabla u\,\nabla \phi dx=\int_{\Omega
}H(x,u,v,\nabla u,\nabla v)\phi \,dx\,\,\,\forall \phi \in
H_{0}^{1}(\Omega ), \\
\mbox{} \\
\displaystyle\int_{\Omega }\nabla v\,\nabla \psi dx=\int_{\Omega
}G(x,u,v,\nabla u,\nabla v)\psi \,dx\,\,\,\forall \psi \in
H_{0}^{1}(\Omega ),%
\end{array}%
\end{equation*}%
that is, $(u,v)$ is a solution of the system
\begin{equation*}
\left\{
\begin{array}{l}
-\Delta {u}=H(x,u,v,\nabla u,\nabla v)\,\,\,\mbox{in}\,\,\,\Omega , \\
\mbox{} \\
-\Delta {v}=G(x,u,v,\nabla u,\nabla v)\,\,\,\mbox{in}\,\,\,\Omega , \\
\mbox{} \\
u(x)=\tilde{g},v(x)=\widehat{g}\,\,\,\mbox{on}\,\,\,\partial \Omega .%
\end{array}%
\right. \leqno{(AS)}
\end{equation*}
\end{theorem}

\noindent \textbf{Proof.} Here, we will adapt some arguments found in Leon \cite{leon}.  Firstly, we introduce two new functions
$$
H_{1}(x,s,t,\eta ,\xi )=
\left\{
\begin{array}{l}
H(x,\underline{u}(x),\underline{v}(x),\nabla \underline{u}(x),\nabla
\underline{v}(x)),\,\,s\leq \underline{u}(x) \\
H(x,s,\underline{v}(x),\eta ,\nabla \underline{v}(x)),\,\,\,\underline{u}%
(x)\leq s\leq \overline{u}(x)\text{ and }t\leq \underline{v}(x) \\
H(x,s,t,\eta ,\xi ),\,\,\,\underline{u}(x)\leq s\leq \overline{u}(x)\text{
and }\underline{v}(x)\leq t\leq \overline{v}(x) \\
H(x,s,\overline{v}(x),\eta ,\nabla \overline{v}(x)),\,\,\,\underline{u}%
(x)\leq s\leq \overline{u}(x)\text{ and }t\geq \overline{v}(x) \\
H(x,\overline{u}(x),\overline{v}(x),\nabla \overline{u}(x),\nabla \overline{u%
}(x)),\,\,s\geq \overline{u}(x).
\end{array}
\right.
$$and
$$
G_1(x,s,t,\eta,\xi)=
\left\{
\begin{array}{l}
G(x,\underline{u}(x),\underline{v}(x), \nabla \underline{u}(x), \nabla \underline{v}(x)), \,\, t \leq \underline{v}(x) \\
G(x,\underline{u}(x) ,t,\nabla \underline{u}(x), \xi), \,\,\, \underline{v}(x) \leq t \leq \overline{v}(x)\text{ and } s \leq \underline{u}(x) \\
G(x,s,t,\eta, \xi), \,\,\, \underline{v}(x) \leq t \leq \overline{v}(x)\text{ and }  \underline{u}(x) \leq s \leq \overline{u}(x)  \\
G(x,\overline{u}(x),t,\nabla \overline{u}(x), \xi), \,\,\, \underline{v}(x) \leq t \leq \overline{v}(x)\text{ and } s \geq \overline{u}(x) \\
G(x,\overline{u}(x),\overline{v}(x), \nabla \overline{u}(x), \nabla \overline{v}(x)), \,\, t \geq \overline{v}(x). \\ 
\end{array}
\right.
$$
Moreover, for each $l \in (0,1)$, we consider
$$
\gamma_1(x,s)=-(( \underline{u}(x)-s)_{+})^{l}+(( s-\overline{u}(x))_{+})^{l}
$$
and
$$
\gamma_2(x,t)=-(( \underline{v}(x)-t)_{+})^{l}+(( t-\overline{v}(x))_{+})^{l}.
$$
Using the above functions, we will work with the ensuing auxiliary system 

\begin{equation*}
\left\{
\begin{array}{l}
-\Delta {u}=H_1(x,u,v,\nabla u,\nabla v)- \gamma_1(x,u) \,\,\,\mbox{in}\,\,\,\Omega , \\
\mbox{} \\
-\Delta {v}=G_1(x,u,v,\nabla u,\nabla v) - \gamma_2(x,v) \,\,\,\mbox{in}\,\,\,\Omega , \\
\mbox{} \\
u(x)=\tilde{g},v(x)=\widehat{g}\,\,\,\mbox{on}\,\,\,\partial \Omega .%
\end{array}%
\right. \leqno{(S_1)}
\end{equation*}
Setting the functions  
$$
H_2(x,s,t,\eta,\xi)=H_1(x,s,t,\eta,\xi)- \gamma_1(x,s)
$$
and
$$
G_2(x,s,t,\eta,\xi)=G_1(x,s,t,\eta,\xi)- \gamma_2(x,t),
$$
we define the operator  $B:E \to E'$ given by
$$
\begin{array}{l}
\left\langle B(u,v),(\phi,\psi)\right\rangle=\displaystyle \int_{\Omega}(\nabla u \nabla \phi + \nabla v \nabla \psi)dx  - \int_{\Omega}H_2(x,u,v,\nabla u,\nabla v)\phi\,dx  \\
\mbox{}\\
\;\;\;\;\;\;\,\,\,\,\,\,\;\;\;\;\;\;\;\;\;\;\;\;\;\;\;\; - \displaystyle \int_{\Omega}G_2(x,u,v,\nabla u,\nabla v)\psi\,dx,
\end{array}
$$
where $E=H^{1}_{0}(\Omega) \times H^{1}_{0}(\Omega)$ is endowed of the norm
$$
\|(u,v)\|=(\|u\|^{2}+\|v\|^{2})^{\frac{1}{2}},
$$
with $\|\mbox{\,\,\,}\|$ being the usual norm in $H^{1}_{0}(\Omega)$. 

Using the hypotheses on $H$ and $G$ together with the definition of $H_1,H_2,G_1, G_2, \gamma_1$ and $\gamma_2$, we can prove the ensuing properties for operator $B$: \\

\noindent {$\bf I) \quad  B$ \textbf{is continuous:}   \\

\noindent The proof of this property follows by using the fact that $H_1,G_1$ belong to $L^{\infty}$. \\

\noindent {$\bf II) \quad  B$ \textbf{is bounded:}   \\

\noindent Here, the boundedness of $B$ is understood in the sense that if $U \subset E$ is a bounded set, then $B(U) \subset E'$ is also bounded. This property also follows using the boundedness of $H_1$ and $G_1$. \\

\noindent {$\bf III) \quad  B$ \textbf{is coercive:}   \\

\noindent Here, it is enough to prove that 
$$
\frac{\left\langle B(u,v),(u,v)  \right\rangle  }{\|(u,v)\|} \to +\infty \quad \mbox{as} \quad \|(u,v)\| \to +\infty .
$$
Using again the boundedness of $H_1$ and $G_1$, we derive 
$$
\left\langle B(u,v),(u,v)  \right\rangle \geq \|(u,v)\|^{2} -C_1 \|(u,v)\|-C_2\|(u,v)\|^{l+1}.
$$
Thus,
$$
\frac{\left\langle B(u,v),(u,v)  \right\rangle  }{\|(u,v)\|} \geq \|(u,v)\| -C_1 - C_2 \|(u,v)\|^{l},
$$
showing that 
$$
\frac{\left\langle B(u,v),(u,v)  \right\rangle  }{\|(u,v)\|} \to +\infty \quad \mbox{as} \quad \|(u,v)\| \to +\infty .
$$

\noindent {$\bf IV) \quad  B$ \textbf{is pseudomonotone:} \\

First of all, we recall that $B$ is a pseudomonotone operator if \linebreak $(u_n,v_n) \rightharpoonup (u,v) $ in $E$ and verifies  
\begin{equation} \label{B1}
\limsup_{n \to +\infty} \left\langle B(u_n,v_n), (u_n,v_n)-(u,v)  \right\rangle \leq 0,
\end{equation}
then 
\begin{equation} \label{B2}
\liminf_{n \to +\infty}\left\langle B(u_n,v_n), (u_n,v_n)-(\phi,\psi)\right\rangle \geq \left\langle B(u,v), (u,v)-(\phi, \psi) \right\rangle \,\, \forall (\phi, \psi) \in E.
\end{equation}
In our case, the weak limit $(u_n,v_n) \rightharpoonup (u,v)$ in $E$ yields 
$$
\int_{\Omega}H_1(x,u_n,v_n,\nabla u_n, \nabla v_n)(u_n-u) \to 0
$$
and
$$
\int_{\Omega}G_1(x,u_n,v_n,\nabla u_n, \nabla v_n)(v_n-v) \to 0.
$$
Thereby, the above limits combined with  (\ref{B1}) load to 
$$
\limsup_{n \to +\infty}\int_{\Omega}(\nabla u_n \nabla (u_n-u) + \nabla v_n \nabla (v_n-v)) \leq 0 ,
$$
from where it follows that 
$$
(u_n,v_n) \to (u,v) \quad \mbox{in} \quad E.
$$

The properties $I)-IV)$ allow us to use \cite[Theorem 3.3.6]{Necas} to conclude that $B$ is surjective. Therefore,  there exists $(u,v) \in E$ such that 
$$
\left\langle B(u,v),(\phi,\psi) \right\rangle=0 \,\,\ \forall (\phi,\psi) \in E,
$$ 
implying that $(u,v)$ is a solution of $(S_1)$. Now, our goal is showing that
\begin{equation} \label{B3}
(i) \,\,\,  \underline{u} \leq u \leq \overline{u} \quad \mbox{and} \quad (ii) \,\,\, \underline{v} \leq v \leq \overline{v}.
\end{equation}
We will show only $(i)$, because the same arguments can be used to prove $(ii)$.   Choosing $(\phi, \psi)=((u-\overline{u})_{+},0)$ as a test function, we have 
$$
\int_{\Omega}\nabla u \nabla (u-\underline{u})_{+}=\int_{\Omega} H_2(x,u,v,\nabla u, \nabla v)(u-\overline{u})_{+}\,dx.
$$ 
From definition of $H_2$,
$$
\int_{\Omega}\nabla u \nabla (u-\overline{u})_{+}=\int_{\Omega}H_1(x,u,v,\nabla u, \nabla v)(u-\overline{u})_{+}\,dx - \int_{\Omega}\gamma_1(x,u)(u-\overline{u})_{+}\,dx
$$ 
and so, 
$$
\int_{\Omega}\nabla u \nabla (u-\overline{u})_{+}dx=\int_{\Omega}H(x,\overline{u},\overline{v},\nabla \overline{u}, \nabla \overline{v})(u-\overline{u})_{+}\,dx - \int_{\Omega}(u-\overline{u})_{+}^{l+1}\,dx.
$$ 
Since $(\overline{u},\overline{v})$ is a supersolution, it follows that 
$$
\int_{\Omega}\nabla u \nabla (u-\overline{u})_{+}dx \leq \int_{\Omega}\nabla \overline{u} \nabla (u-\overline{u})_{+}\,dx - \int_{\Omega}(u-\overline{u})_{+}^{l+1}\,dx,
$$ 
or equivalently,
$$
\int_{\Omega}|\nabla (u-\overline{u})_{+}|^{2}dx \leq - \int_{\Omega}(u-\overline{u})_{+}^{l+1}\,dx \leq 0,
$$
showing that $(u-\overline{u})_{+}=0$, from where it follows that $u \leq \overline{u}$.  To prove that $\underline{u} \leq u$, we choose $(\phi, \psi)=((\underline{u}-u)_{+},0)$ as a test function. Repeating the above arguments, we get 
$$
\int_{\Omega}\nabla u \nabla (\underline{u}-u)_{+}dx=\int_{\Omega}H(x,\underline{u},\underline{v},\nabla \underline{u}, \nabla \underline{v})(\underline{u}-u)_{+}\,dx +\int_{\Omega}(\underline{u}-u)_{+}^{l+1}\,dx.
$$ 
Since $(\underline{u},\underline{v})$ is a subsolution, it follows that 
$$
\int_{\Omega}\nabla u \nabla (\underline{u}-u)_{+} dx \geq \int_{\Omega}\nabla \underline{u} \nabla (\underline{u}-u)_{+}\,dx + \int_{\Omega}(\underline{u}-u)_{+}^{l+1}\,dx,
$$ 
or equivalently,
$$
\int_{\Omega}|\nabla (\underline{u}-u)_{+}|^{2} dx \leq - \int_{\Omega}(\underline{u}-u)_{+}^{l+1}\,dx \leq 0,
$$
showing that $(\underline{u}-u)_{+}=0$, and so, $\underline{u} \leq u$. 
Combining (\ref{B3}) with the definition of $H_2$ and $G_2$, it follows that
$$
H_2(x,u,v,\nabla u, \nabla v)= H(x,u,v,\nabla u, \nabla v)
$$
and 
$$
G_2(x,u,v,\nabla u, \nabla v)= G(x,u,v,\nabla u, \nabla v),
$$
showing that $(u,v)$ is a solution for system $(AS)$. \hfill\rule{2mm}{2mm}

\Section{Existence of solution for system $(S)_{-}$}

\label{S3}

In this section, we will study the existence of solution for the following
singular elliptic system
\begin{equation*}
\left\{
\begin{array}{l}
-\Delta{u} = \frac{1}{v^{\alpha_1}}-v^{\beta_1}+ g_1(\nabla u, \nabla v)
\,\,\, \mbox{in} \,\,\, \Omega, \\
\mbox{} \\
-\Delta{v} =\frac{1}{u^{\alpha_2}}-u^{\beta_2}+ g_2(\nabla u, \nabla v)
\,\,\, \mbox{in} \,\,\, \Omega, \\
\mbox{} \\
u,v>0 \,\,\ \mbox{in} \,\,\, \Omega, \\
\mbox{} \\
u=v=0 \,\,\  \mbox{on} \,\,\, \partial   \Omega , \\
\end{array}
\right. \leqno{(S)_{-}}
\end{equation*}
where $\Omega$ is a bounded domain in $\mathbb{R}^{N}$ with smooth boundary,
\linebreak $g_i: \mathbb{R}^{2N} \to [0, +\infty)$ are positive continuous
functions belong to $L^{\infty}(\mathbb{R}^{2N})$ and $\alpha_i,\beta_i \in
[0,1).$

Our approach consists in considering for $\epsilon >0$ the approximated
system
\begin{equation*}
\left\{
\begin{array}{l}
-\Delta {u}=\frac{1}{(|v|^{2}+\epsilon )^{\frac{\alpha _{1}}{2}}}-v^{\beta
_{1}}+g_{1}(\nabla u,\nabla v)\,\,\,\mbox{in}\,\,\,\Omega , \\
\mbox{} \\
-\Delta {v}=\frac{1}{(|u|^{2}+\epsilon )^{\frac{\alpha _{2}}{2}}}-u^{\beta
_{2}}+g_{2}(\nabla u,\nabla v)\,\,\,\mbox{in}\,\,\,\Omega , \\
\mbox{} \\
u=v=0\,\,\  \mbox{on} \,\,\,  \partial  \Omega .%
\end{array}%
\right. \leqno{(S_\epsilon)_{-}}
\end{equation*}

For this class of system it is possible to find sub and supersolution which
do not depend of $\epsilon $. For example, by using the positive
eigenfunction associated with the first eigenvalue, we can find easily a
subsolution $(\underline{u},\underline{u})$. To get the supersolution, we
observe that any large constant $M>0$ can be used to get a supersolution $(%
\overline{u},\overline{v})=(M,M)$ with $M>\Vert \underline{u}\Vert _{\infty
} $.

Now, let us consider the functions 
$$
H_{\epsilon}(x,s,t,\eta,\xi)= \frac{1}{(|t|^{2}+\epsilon )^{\frac{\alpha _{1}}{2}}}-t^{\beta
_{1}}+g_{1}(\eta,\xi)
$$
and
$$
G_{\epsilon}(x,s,t,\eta,\xi)=\frac{1}{(|s|^{2}+\epsilon )^{\frac{\alpha _{2}}{2}}}-s^{\beta
_{2}}+g_{2}(\eta,\xi),
$$
which are well defined in $\Omega \times \mathbb{R}^{+}\times \mathbb{R}^{+}\times \mathbb{R}%
^{N}\times \mathbb{R}^{N}$ . 

By using Theorem \ref{T2}, there exists a positive $(u_{\epsilon
},v_{\epsilon })$ verifying system $(S_{\epsilon })_{-}$. Set $\epsilon =%
\frac{1}{n}$ with any integer $n\geq 1$. From now on, we denote by $%
(u_{n},v_{n})$ the solution $(u_{\frac{1}{n}},v_{\frac{1}{n}})$. Hence,%
\begin{equation*}
\left\{
\begin{array}{l}
-\Delta {u_{n}}=\frac{1}{(|v_{n}|^{2}+\frac{1}{n})^{\frac{\alpha _{1}}{2}}}%
-v_{n}^{\beta _{1}}+g_{1}(\nabla u_{n},\nabla v_{n})\,\,\,\mbox{in}%
\,\,\,\Omega , \\
\mbox{} \\
-\Delta {v_{n}}=\frac{1}{(|u_{n}|^{2}+\frac{1}{n})^{\frac{\alpha _{2}}{2}}}%
-u_{n}^{\beta _{2}}+g_{2}(\nabla u_{n},\nabla v_{n})\,\,\,\mbox{in}%
\,\,\,\Omega , \\
\mbox{} \\
u_{n}=v_{n}=0\,\, \mbox{on} \,\,\, \partial   \Omega .%
\end{array}%
\right. \leqno{(S_n)_{-}}
\end{equation*}

Once that $\alpha _{i},\beta _{i}\in \lbrack 0,1)$, $i=1,2$, and $g_{i}$
belongs to $L^{\infty }(\mathbb{R}^{2N})$, it follows that $(u_{n},v_{n})$
is bounded in $H_{0}^{1}(\Omega )\times H_{0}^{1}(\Omega )$. Indeed, since
\begin{equation*}
\begin{array}{c}
u_{n}\geq \underline{u}>0\text{ and }v_{n}\geq \underline{v}>0\text{\ in }%
\Omega ,%
\end{array}%
\end{equation*}%
we have%
\begin{equation*}
\Vert u_{n}\Vert ^{2}\leq \int_{\Omega }\frac{u_{n}}{v_{n}^{\alpha _{1}}}%
dx+\int_{\Omega }g_{1}(\nabla u_{n},\nabla v_{n})u_{n}\text{ }dx
\end{equation*}%
and so,
\begin{equation}
\Vert u_{n}\Vert ^{2}\leq \int_{\Omega }\frac{u_{n}}{\underline{v}^{\alpha
_{1}}}dx+\Vert g_{1}\Vert _{\infty }\int_{\Omega }|u_{n}|\,dx.  \label{4}
\end{equation}%
Similarly we derive
\begin{equation}
\displaystyle\left\Vert v_{n}\right\Vert ^{2}\leq \int_{\Omega }\frac{v_{n}}{%
\underline{u}^{\alpha _{2}}}dx+\Vert g_{2}\Vert _{\infty }\int_{\Omega
}|v_{n}|\,dx.  \label{5}
\end{equation}%
On the other hand, for $\alpha _{i}\in \lbrack 0,1)$, $i=1,2$, we may invoke
the Hardy-Sobolev inequality in the form stated in \cite[Lemma 2.3]%
{AlvesCorrea} to infer that
\begin{equation}
\int_{\Omega }\frac{u_{n}}{\underline{v}^{\alpha _{1}}}\,dx\leq C \Vert
u_{n}\Vert \,\,\,\mbox{and}\,\,\,\int_{\Omega }\frac{v_{n}}{\underline{u}%
^{\alpha _{2}}}\,dx\leq C \Vert v_{n}\Vert .  \label{6}
\end{equation}%
Combining (\ref{4})-(\ref{6}) with Sobolev embedding, it follows that $%
(u_{n},v_{n})$ is bounded in $H_{0}^{1}(\Omega )\times H_{0}^{1}(\Omega )$.
Consequently, we can assume that there is $(u,v)\in H_{0}^{1}(\Omega )\times
H_{0}^{1}(\Omega )$ and $G_{i}\in L^{2}(\Omega )$ verifying
\begin{equation}
u_{n}\rightharpoonup u,v_{n}\rightharpoonup v\,\,\,\mbox{in}%
\,\,\,H_{0}^{1}(\Omega ),  \label{1}
\end{equation}%
\begin{equation}
u_{n}\rightarrow u,v_{n}\rightarrow v\,\,\,\mbox{in}\,\,\,L^{p}(\Omega
)\,\,\,\text{for all }p\in \lbrack 1,+\infty ),  \label{1.1}
\end{equation}%
\begin{equation}
u_{n}(x)\rightarrow u(x),v_{n}(x)\rightarrow v(x)\,\,\,\mbox{a.e. in}%
\,\,\,\Omega  \label{1.2}
\end{equation}%
and
\begin{equation}
g_{i}(\nabla u_{n},\nabla v_{n})\rightharpoonup G_{i}\,\,\,\mbox{in}%
\,\,\,L^{2}(\Omega ).  \label{2}
\end{equation}%
Once that
\begin{equation*}
\begin{array}{c}
\underline{u}\leq u_{n}\leq M\text{ and }\underline{v}\leq v_{n}\leq M\text{
a.e in }\Omega ,\,\,\text{for all }n\in \mathbb{N},%
\end{array}%
\end{equation*}%
the limit (\ref{1.2}) gives
\begin{equation*}
\begin{array}{c}
\underline{u}\leq u\leq M\text{ and }\underline{v}\leq v\leq M\,\,\text{a.e
in }\Omega .%
\end{array}%
\end{equation*}%
Recall that $(S_{n})_{-}$ entails%
\begin{equation}
\left\{
\begin{array}{l}
\displaystyle\int_{\Omega }\nabla u_{n}\nabla \varphi \,dx=\int_{\Omega }\Big(
\frac{1}{(|v_{n}|^{2}+\frac{1}{n})^{\frac{\alpha _{1}}{2}}}-v_{n}^{\beta
_{1}}+g_{1}(\nabla u_{n},\nabla v_{n})\Big)\varphi \ dx \\
\displaystyle\int_{\Omega }\nabla u_{n}\nabla \psi \,dx=\int_{\Omega }\Big(\frac{%
1}{(|u_{n}|^{2}+\frac{1}{n})^{\frac{\alpha _{2}}{2}}}-u_{n}^{\beta
_{2}}+g_{2}(\nabla u_{n},\nabla v_{n})\Big)\psi \ dx%
\end{array}%
\right.  \label{3}
\end{equation}%
for all $(\varphi ,\psi )\in H_{0}^{1}(\Omega )\times H_{0}^{1}(\Omega ).$
Setting $(\varphi ,\psi )=(u_{n}-u,v_{n}-v)$ in (\ref{3}) yields%
\begin{equation*}
\left\{
\begin{array}{c}
\displaystyle\int_{\Omega }\nabla u_{n}\nabla (u_{n}-u)\,dx=\int_{\Omega }\Big(%
\frac{1}{(|v_{n}|^{2}+\frac{1}{n})^{\frac{\alpha _{1}}{2}}}-v_{n}^{\beta
_{1}}+g_{1}(\nabla u_{n},\nabla v_{n})\Big)(u_{n}-u)\ dx \\
\displaystyle\int_{\Omega }\nabla v_{n}\nabla (v_{n}-v)\,dx=\int_{\Omega }\Big(%
\frac{1}{(|u_{n}|^{2}+\frac{1}{n})^{\frac{\alpha _{2}}{2}}}-u_{n}^{\beta
_{2}}+g_{2}(\nabla u_{n},\nabla v_{n})\Big)\text{ }(v_{n}-v)\ dx.%
\end{array}%
\right.
\end{equation*}%
We point out that (\ref{1})-(\ref{2}) implies that
\begin{equation*}
\lim_{n\rightarrow \infty }\int_{\Omega }\Big(\frac{1}{(|v_{n}|^{2}+\frac{1}{n}%
)^{\frac{\alpha _{1}}{2}}}-v_{n}^{\beta _{1}}+g_{1}(\nabla u_{n},\nabla
v_{n}\Big))(u_{n}-u)\ dx=0
\end{equation*}%
and
\begin{equation*}
\lim_{n\rightarrow \infty }\int_{\Omega }\Big(\frac{1}{(|u_{n}|^{2}+\frac{1}{n}%
)^{\frac{\alpha _{2}}{2}}}-u_{n}^{\beta _{2}}+g_{2}(\nabla u_{n},\nabla
v_{n})\Big)\text{ }(v_{n}-v)\ dx=0.
\end{equation*}%
Consequently,
\begin{equation*}
\begin{array}{c}
\displaystyle\lim_{n\rightarrow +\infty }\int_{\Omega }\nabla u_{n}\nabla
(u_{n}-u)\,dx=\lim_{n\rightarrow +\infty }\int_{\Omega }\nabla v_{n}\nabla
(v_{n}-v)\,dx=0%
\end{array}%
\end{equation*}%
leading to
\begin{equation}
(u_{n},v_{n})\rightarrow (u,v)\text{ in }H_{0}^{1}(\Omega )\times
H_{0}^{1}(\Omega ).  \label{13}
\end{equation}%
This way, passing to relabeled subsequences, we have the limits
\begin{equation*}
\nabla u_{n}(x)\rightarrow \nabla u(x)\,\,\,\mbox{and}\,\,\,\nabla
v_{n}(x)\rightarrow \nabla v(x)\,\,\,\mbox{a.e in}\,\,\,\Omega ,
\end{equation*}%
which imply that
\begin{equation}
g_{i}(\nabla u_{n},\nabla v_{n})\rightharpoonup g_{i}(\nabla u,\nabla
v)\,\,\,\mbox{in}\,\,\,L^{2}(\Omega )\,\,(i=1,2).  \label{14}
\end{equation}%
Now, from (\ref{13}) and (\ref{14}), we may pass to the limit in (\ref{3}) to
conclude that $(u,v)$ is a solution for $(S)_{-}$. This completes the proof.

\Section{Existence of solution for system $(S)_{+}$}

\label{S4}

In this section, we will study the existence of solution for the following
singular elliptic system
\begin{equation*}
\left\{
\begin{array}{l}
-\Delta {u}=\frac{1}{v^{\alpha _{1}}}+v^{\beta _{1}}+g_{1}(\nabla u,\nabla
v)\,\,\,\mbox{in}\,\,\,\Omega , \\
\mbox{} \\
-\Delta {v}=\frac{1}{u^{\alpha _{2}}}+u^{\beta _{2}}+g_{2}(\nabla u,\nabla
v)\,\,\,\mbox{in}\,\,\,\Omega \\
\mbox{} \\
u,v>0\,\,\ \mbox{in}\,\,\,\Omega , \\
\mbox{} \\
u=v=0\,\,\ \mbox{on} \,\,\, \partial   \Omega , \\
\end{array}%
\right. \leqno{(S)_{+}}
\end{equation*}%
by supposing the same hypotheses of Section \ref{S3}.

The existence of solution for $(S)_{+}$ can be obtained by using the same
arguments explored in the previous section. The unique difference is in the
construction of the supersolution that we will work. Here, the idea is the
following:

Fix $R>0$ large enough such that $\Omega \subset B_{R}(0)$ and denote by $e$
the unique solution of the problem
\begin{equation*}
\left\{
\begin{array}{l}
-\Delta e=1,\,\,\,\mbox{in}\,\,\,B_{R}(0) \\
\mbox{} \\
e=0,\,\,\mbox{on}\,\,\partial B_{R}(0).%
\end{array}%
\right.
\end{equation*}%
Recalling that $g_{i}\in L^{\infty }$, for $M>\Vert \underline{u}\Vert
_{\infty }$ large enough, a simple computation shows that
\begin{equation*}
\left\{
\begin{array}{l}
-\Delta (Me)=M\geq \frac{1}{(Me)^{\alpha _{1}}}+(Me)^{\beta
_{1}}+g_{1}(\nabla w_{1},\nabla w_{2})\,\,\,\mbox{in}\,\,\,\Omega , \\
\mbox{} \\
-\Delta {(Me)}=M\geq \frac{1}{(Me)^{\alpha _{2}}}+(Me)^{\beta
_{2}}+g_{2}(\nabla w_{1},\nabla w_{2})\,\,\,\mbox{in}\,\,\,\Omega , \\
\mbox{} \\
(Me)>0\,\,\ \mbox{in}\,\,\,\overline{\Omega },%
\end{array}%
\right. \leqno{(S)_{+}}
\end{equation*}%
for any $w_1,w_2 \in H^{1}_{0}(\Omega)$. Thereby, the pairs $(\underline{u},%
\underline{u})$ and $(\overline{u},\overline{v})=(Me,Me)$ satisfy the
hypotheses of Theorem \ref{T1}. \hfill \rule{2mm}{2mm}

\vspace{0.5 cm}

\noindent \textbf{Acknowledgement.} This work was accomplished while the second author
was visiting Laboratoire de Math\'{e}matiques, Physique et Syst\`emes
(LAMPS) of Perpignan University, with Averro\`es fellowship. He wishes to
thank Perpignan University for the kind hospitality.

\end{document}